\documentclass[12pt,a4paper,english,reqno]{amsart}
\usepackage[a4paper,footskip=1.5em]{geometry}
\usepackage{amsmath,amssymb,amsthm,mathtools,mathrsfs}
\usepackage{adjustbox,tikz,calc,graphics,babel}
\usepackage{csquotes,enumerate,verbatim}
\usepackage[final]{microtype}
\usepackage[numbers]{natbib}
\usetikzlibrary{shapes.misc,calc,intersections,patterns,decorations.pathreplacing}
\usepackage{hyperref}
\hypersetup{colorlinks=true,linkcolor=blue,citecolor=blue}

\theoremstyle{plain}
\newtheorem*{theorem*}{Theorem}
\newtheorem{theorem}{Theorem}[section]
\newtheorem{lemma}[theorem]{Lemma}

\newtheorem{proposition}[theorem]{Proposition}
\newtheorem*{claim*}{Claim}

\newtheorem{problem}[theorem]{Problem}
\theoremstyle{remark}
\newtheorem*{remark}{Remark}

\def\N{\mathbb{N}}

\def\P{\mathbb{P}}
\def\E{\mathbb{E}}
\def\C{\mathcal}

\def\Scr{\mathscr}

\DeclareMathOperator\Deg{d}
\DeclareMathOperator\Bi{Bin}

\let\emptyset\varnothing
\let\eps\varepsilon

\let\originalleft\left
\let\originalright\right
\renewcommand{\left}{\mathopen{}\mathclose\bgroup\originalleft}
\renewcommand{\right}{\aftergroup\egroup\originalright}

\makeatletter
\def\imod#1{\allowbreak\mkern10mu({\operator@font mod}\,\,#1)}
\makeatother

\begin{document}

\title{Disjoint induced subgraphs of the same order and size}

\author{B\'{e}la Bollob\'{a}s}
\address{Department of Pure Mathematics and Mathematical Statistics, University of Cambridge, Wilberforce Road, Cambridge CB3\thinspace0WB, UK; {\em and\/}
Department of Mathematical Sciences, University of Memphis, Memphis TN 38152, USA; {\em and\/} London Institute for Mathematical Sciences, 35a South St., Mayfair, London W1K\thinspace2XF, UK.}
\email{b.bollobas@dpmms.cam.ac.uk}

\author{Teeradej Kittipassorn}
\address{Department of Mathematical Sciences, University of Memphis, Memphis TN 38152, USA}
\email{t.kittipassorn@memphis.edu}

\author{Bhargav Narayanan}
\address{Department of Pure Mathematics and Mathematical Statistics, University of Cambridge, Wilberforce Road, Cambridge CB3\thinspace0WB, UK}
\email{b.p.narayanan@dpmms.cam.ac.uk}

\author{Alex Scott}
\address{Mathematical Institute, University of Oxford, Andrew Wiles Building, Radcliffe Observatory Quarter, Woodstock Road, Oxford OX2\thinspace6GG, UK}
\email{scott@maths.ox.ac.uk}

\date{5 January 2014}

\begin{abstract}
For a graph $G$, let $f(G)$ be the largest integer $k$ for which there exist two vertex-disjoint induced subgraphs of $G$ each on $k$ vertices, both inducing the same number of edges. We prove that $f(G) \ge n/2 - o(n)$ for every graph $G$ on $n$ vertices. This answers a question of Caro and Yuster.
\end{abstract}

\maketitle

\section{Introduction}
Given a graph $G$, can we guarantee that $G$ contains two large, vertex-disjoint copies of the same graph? It follows from Ramsey's theorem that any graph on $n$ vertices contains two vertex-disjoint isomorphic induced subgraphs on $\Omega( \log{n})$ vertices; by considering a random graph on $n$ vertices, it is easy to check that this is also best-possible up to constant factors.

What if, rather than asking for two isomorphic subgraphs, we ask for two subgraphs that are the same with respect to one or more graph parameters? Caro and Yuster~\citep{CarYus} considered the question of finding two vertex-disjoint subgraphs of a given graph of the same order which induce the same number of edges. For a graph $G$, let $f(G)$ be the largest integer $k$ such that there are two vertex-disjoint induced subgraphs of $G$ each on $k$ vertices, both inducing the same number of edges and let $f(n)$ be the minimum value of $f(G)$ taken over all graphs on $n$ vertices. Trivially, $f(n) \le \lfloor n/2 \rfloor$; also, as shown by Ben-Eliezer and Krivelevich~\citep{BEKrivel}, equality holds (with high probability) for the Erd\H{o}s--R\'enyi random graphs $G(n,p)$ for all $0\le p \le 1$.

There was a large gap between the best known upper and lower bounds for $f(n)$. From below, one can easily show using the pigeonhole principle that $f(n) = \Omega (n^{1/3})$. As observed by Caro and Yuster, it is possible to improve this to $f(n) = \Omega (n^{1/2})$ using a well known result of Lov\'asz determining the chromatic number of Kneser graphs. By considering a carefully constructed disjoint union of cliques, each on an odd number of vertices, Caro and Yuster showed that $f(n) \le n/2 - \Omega(\log \log n)$.

As expected, one can say more about $f(G)$ when $G$ belongs to certain special graph classes. For example, Axenovich, Martin and Ueckerdt~\citep{AxeMarU} showed that $f(G) \ge \lceil n/2 \rceil-1$ when $G$ is a forest; this is clearly best-possible. Indeed, it is possible to get quite close to the trivial upper bound of $n/2$ when we restrict our attention to sparse graphs. In their paper, Caro and Yuster showed, for any fixed $\alpha > 0$, that if $G$ is a graph on $n$ vertices, then $f(G) \ge n/2 -o(n)$ provided $G$ has at most $n^{2-\alpha}$ edges (or non-edges).  Axenovich, Martin and Ueckerdt~\citep{AxeMarU} later showed that the same holds for graphs with at most $o(n^2/ (\log n)^2)$ edges.

Our main aim in this paper is to narrow considerably the gap between the best known upper and lower bounds for $f(n)$, and thereby answer a question of Caro and Yuster~\citep{CarYus}.

\begin{theorem}
\label{subedge13-mainthm}
For every $ \eps >0$, there exists a natural number $N=N(\eps)$ such that for any graph $G$ on $n>N$ vertices, $f(G) \ge n/2 - \eps n$. Consequently, 
\[n/2 - o(n) \le f(n) \le n/2 - \Omega(\log \log n).\]
\end{theorem}

We remark that much research has been done on the family of induced subgraphs of a graph. For example, call a graph {\em $k$-universal} if it contains every graph of order $k$ as an induced subgraph. Very crudely, if $G$ is a $k$-universal graph with $n$ vertices, then 
\[\binom{n}{k}\ge \frac{2^{\binom{k}{2}}}{k!},\] 
so $n\ge 2^{(k-1)/2}$. As remarked in~\citep{BT}, almost all graphs with $k^2 2^{k/2}$ vertices are $k$-universal, and the Paley graphs come close to providing examples which are almost as good. Hajnal conjectured that if a graph only has a `small' number of distinct (non-isomorphic) induced subgraphs, then it contains a trivial (complete or empty) subgraph with linearly many vertices. This was proved, shortly after the conjecture was made, by Alon and Bollob\'{a}s~\citep{AlBo}, and Erd\H{o}s and Hajnal~\citep{ErdHaj1}, the latter in a stronger form. In~\citep{AlBo} only a few parameters, like order, size and maximal degree, were used to distinguish non-isomorphic graphs. 

Erd\H{o}s and Hajnal~\cite{EH-conj} then went much further: they realised that forbidding a single graph as an induced subgraph severely constrains the structure of a graph. More precisely, they made the major conjecture that for every graph $H$, there is a positive constant $\gamma(H)$ such that if a graph of order $n$ does not contain $H$ as an induced subgraph, then the graph contains a trivial subgraph with at least $n^{\gamma(H)}$ vertices. In spite of all the work on this conjecture, see ~\citep{ChuSaf,FoxSud,ProRod} for instance, we are very far from the desired bound. 

Let us finally mention another interesting line of research about finding disjoint isomorphic (not necessarily induced) subgraphs. Jacobson and Sch\"{o}nheim (see~\citep{ErdPachPyb, LeeLohSud}) independently raised the question of finding edge-disjoint isomorphic subgraphs.  Improving on results of Erd\H{o}s, Pyber and Pach~\citep{ErdPachPyb}, it has been shown by Lee, Loh and Sudakov~\citep{LeeLohSud} that every graph on $m$ edges contains a pair of edge-disjoint isomorphic subgraphs with at least $\Omega((m \log{m})^{2/3})$ edges and that this is also best-possible up to a multiplicative constant.

The rest of this paper is organised as follows. We give an overview of our approach in Section~\ref{subedge13-strat}, and then fill in the details and prove Theorem~\ref{subedge13-mainthm} in Section~\ref{subedge13-proof}. There are many natural questions about induced subgraphs which are close to Theorem~\ref{subedge13-mainthm} in spirit; we conclude in Section~\ref{subedge13-conc} by  mentioning some of these.  

\section{Preliminaries}
Our objective in this section is to establish some notational conveniences and collect together, for easy reference, some simple propositions that we shall make use of when proving our main result.

\subsection{Notation}
It will be convenient to establish some notation for working with sets of pairs. A \emph{pair $\{x,y\}$} will always mean an unordered pair with $x \neq y$, and a \emph{collection} of pairs $\C{P}$ will always mean a set of disjoint pairs; for example, $\C{P} = \{ \{1,2\}, \{3,4\}\}$ is a collection of pairs, but $\C{Q} = \{\{1,2\},\{2,3\}\}$ is not. For a collection of pairs denoted by $\mathcal{P}$, we shall write $P$ for the underlying ground set of elements, i.e., $P = \bigcup_{\{x,y\} \in \C{P}} \lbrace x, y \rbrace$; in other words, we reserve the corresponding upper case letter for the ground set. We shall say that two collections of pairs $\C{P}$ and $\C{Q}$ are disjoint if $P \cap Q = \emptyset$; for example, the collections $\C{P}_1 = \{ \{1,2\}, \{3,4\}\}$ and $\C{Q}_1 = \{\{5,6\},\{7,8\}\}$ are disjoint, while the collections $\C{P}_2 = \{ \{1,2\}, \{3,4\}\}$ and $\C{Q}_2 = \{\{1,3\},\{2,4\}\}$ are not.

As usual, given a graph $G=(V,E)$, we write $\Deg(v)$ and $\Gamma (v)$ respectively for the degree and for the neighbourhood of a vertex $v$ in $G$. For a subset $U \subset V$, we write $G[U]$ for the subgraph induced by $U$, $e(U)$ for the number of edges of $G[U]$, and $\Deg(U)$ for the sum of the degrees (in $G$) of the vertices of $U$. Given two disjoint subsets $A,B \subset V$, we write $e(A,B)$ for the number of edges with one endpoint each in $A$ and $B$.

We shall also use the following less common terminology and notation. For any  two vertices $x,y \in V$, we write $\delta(x,y)$ for the \emph{degree difference} between $x$ and $y$, namely the quantity $|\Deg(x)-\Deg(y)|$. We say that two vertices $x$ and $y$  \emph{disagree on a vertex} $v\neq x,y$ if $v$ is adjacent to exactly one of $x$ and $y$; otherwise $x$ and $y$ \emph{agree} on $v$. For any two vertices $x,y \in V$, the \emph{difference neighbourhood} $\Gamma (x,y)$ of $x$ and $y$ is the set of vertices $v\neq x,y$ on which $x$ and $y$ disagree; we write $\Delta (x,y)$ for the size of the difference neighbourhood, so that $\delta(x,y) \le \Delta(x,y)$. If two vertices $x$ and $y$  agree on every vertex $v \neq x,y$, we say that the pair $\{x,y\}$ is a \emph{clone pair}. When the graph $G$ in question is not clear from the context, we shall, for example, write $\delta(x,y,G)$ to denote the degree difference between $x$ and $y$ in $G$.

We say that a graph $G$ is \emph{splittable} if there is a partition $V=A \cup B$ of its vertex set into two sets $A$ and $B$ of equal size with $e(A) = e(B)$; in this case, we call $(A,B)$ a \emph{splitting} of $G$. Note that $e(A) = e(B)$ if and only if $\Deg(A) = \Deg(B)$, since $\Deg(A) = 2e(A)+ e(A,B)$.

Our conventions for asymptotic notation are largely standard; however, we feel obliged to point out that we write $o_{k \to \infty} (1)$ to denote a function (of $k$) that goes to $0$ as $k \to \infty$, and that when we write, say $\Omega_k(.)$, we mean that the constant suppressed by the asymptotic notation is allowed to depend on (but is completely determined by) the parameter $k$. For the sake of clarity of presentation, we systematically omit floors and ceilings whenever they are not crucial.

\subsection{Preliminary observations}
We shall make use of the following simple observation repeatedly when constructing a splitting.

\begin{proposition}
\label{subedge13-approx}
Given positive real numbers $x_1,x_2,\dots,x_t$ in the interval $[a,b]$ with $0\le a\le b$, we may, for every $y \in [-ta,ta]$, choose signs $\zeta_i \in \lbrace -1,+1\rbrace$ such that $| y + \sum \zeta_i x_i | \le  b$. \qed
\end{proposition}

The following first moment bound will prove useful; it is easily checked that the bound is the best-possible.

\begin{proposition}
\label{subedge13-spacesplit}
Let $X$ be a random variable such that $X\le N$ and $\E[X]\ge Np$. Then
\[ \P\left(X\ge \frac{\E[X]}{2}\right) \ge \frac{p}{2-p}.\eqno\qed\]
\end{proposition}

We will also need the following two easy propositions.

\begin{proposition}
\label{subedge13-grouping}
Given $x_1, x_2,\dots,x_t$ in the interval $[0,a]$, a positive real $b$ and a natural number $N$, it is possible to find $\lfloor t/N \rfloor - \lceil a/b \rceil$ disjoint subsets of $\lbrace x_1, x_2,\dots,x_t \rbrace$, each of size $N$, such that $|x_i - x_j|\le b$ for any $x_i$ and $x_j$ belonging to the same subset.
\end{proposition}
\begin{proof}
Suppose that $x_1 \le x_2 \le \dots \le x_t$. Let $i_0=1$ and define $i_j$ to be the smallest index such that $x_{i_j} > x_{i_{j-1}} + b$ and consider the sets $S_j = \lbrace x_{i_j},x_{i_j + 1},\dots,x_{i_{j + 1}-1}\rbrace$. Since $x_1 \ge 0 $ and $x_t \le a$, there are at most $\lceil a/b \rceil$ such sets. Now, by discarding at most $N$ numbers from each $S_j$ if necessary, we can assume that $N$ divides $|S_j|$ for each $j$. We now partition each $S_j$ into subsets of size $N$. Clearly, $|x_i - x_j|\le b$ for any $x_i$ and $x_j$ belonging to the same subset.
The number of elements we have discarded is at most $N \lceil a/b \rceil$.  So the number of subsets of size $N$ we are left with is at least $\lfloor t/N \rfloor - \lceil a/b \rceil$.
\end{proof}

\begin{remark}
We shall often apply Proposition~\ref{subedge13-grouping} to the degrees of a subset of vertices of a graph; we consequently obtain disjoint groups of vertices such that the degree difference of any two vertices in the same group is suitably bounded.
\end{remark}

\begin{proposition}
\label{subedge13-triples}
Let $x$, $y$ and $z$ be three vertices and $U$ some subset of vertices of a graph $G$. Then some two of the vertices $x$, $y$ and $z$ disagree on at most two thirds of the vertices of $U$.
\end{proposition}
\begin{proof}
Any vertex $v \in U$ belongs to at most two of the three difference neighbourhoods $\Gamma (x,y)$, $\Gamma (y,z)$ and $\Gamma (z,x)$. The claim follows by averaging.
\end{proof}

\subsection{Binomial random variables}

We will need some easily proven statements about binomial random variables. We collect these here. As usual, for a random variable with distribution $\Bi (N,p)$, we write $\mu (= Np)$ for its mean and $\sigma^2( = Np(1-p))$ for its variance.

The first proposition we shall require is an easy consequence of the fact that $e^{-2x} \le 1 - x \le e^{-x}$ for all $0 \le x \le 1/2$.

\begin{proposition}
\label{subedge13-binsmall}
Let $X$ be a random variable with distribution $\Bi (N,p)$, with $p\le 1/2$. Then for every $1 \le k\le n$,
\[ \exp \left(-2 \mu \right) \left(\mu/k\right)^k \le \P \left(X = k\right) \le \exp \left(- \mu \right) \left(2e\mu/k\right)^k.\]
Also, $ \exp (-2 \mu ) \le \P (X = 0) \le \exp (- \mu )$. \qed
\end{proposition}

We shall make use of  the following standard concentration result which first appeared in a paper of Bernstein and was later rediscovered by Chernoff and Hoeffding; see~\citep{textbook} for example.

\begin{proposition}
\label{subedge13-bernstein}
Let $X$ be a random variable with distribution $\Bi(N,p)$. Then
\[\P( |X - Np| > t ) \le 2 \exp \left(\frac{-2t^2}{N} \right).\eqno\qed\]
\end{proposition}

\begin{proposition}
\label{subedge13-binparity}
Let $X$ be a random variable with distribution $\Bi(N,p)$. Then
\[\P(X \,\mbox{is even}) = \frac{1}{2}(1+(1-2p)^N).\eqno\qed\]
\end{proposition}

\begin{proposition}
\label{subedge13-binequal}
Let $X_1$ and $X_2$ be two independent random variables both with distribution $\Bi(N,p)$. Then
\[\P(X_1 = X_2) = o_{\sigma \to \infty}(1).\] 
In particular, when $p\le 1/2$, $\P(X_1 = X_2) = o_{\mu \to \infty}(1)$.\qed
\end{proposition}

\begin{proposition}
\label{subedge13-binclose}
Let $X_1$ and $X_2$ be two independent random variables with distributions $\Bi(N_1,p)$ and $\Bi(N_2,p)$ respectively, with $p\ge 1/2$. Then
\[\P(|X_1 - X_2| < |N_1 - N_2|^{1/3} ) = o_{|N_1 - N_2| \to \infty}(1).\eqno\qed\]
\end{proposition}

\begin{proposition}
\label{subedge13-binfar}
Let $X_1$ and $X_2$ be two independent random variables with distributions $\Bi(N_1,p)$ and $\Bi(N_2,p)$ respectively, with $p\ge 1/2$. Suppose $N_1\le N$, $ N_2 \le N$ and $|N_1 - N_2|\le cN^{1/2}$ for some absolute constant $c$. Then
\[\P(|X_1 - X_2| > N^{{2/3}} ) = O \left( \exp \left(  \frac{-N^{1/3}}{5} \right)\right).\eqno\qed\]
\end{proposition}

\section{Overview of our strategy}
\label{subedge13-strat}
To prove Theorem~\ref{subedge13-mainthm}, we need to show that if $\eps > 0$ and $n$ is sufficiently large, then any graph $G$ on $n$ vertices contains two disjoint subsets of vertices of the same size, each of cardinality at least $(1/2 - \eps) n$, which induce the same number of edges. Equivalently, we need to show that it is possible to transform $G$ into a splittable graph by deleting at most $2\eps n$ vertices from $G$. Recall that a graph is splittable if and only if there is a partition of its vertex set into two sets of equal size such that the sums of the degrees of the vertices in the two sets are equal.

We shall show that there is a probability $0 < p \le \eps$ (depending on $G$) such that if we delete vertices from $G$ with probability $p$, then the resulting graph $H$ is splittable with positive probability. 

To show that this random subgraph $H$ is splittable, we shall exhibit a large collection of `gadgets' in $H$. Given $0\le a\le b$, by an \emph{$[a,b]$-gadget}, we mean a pair of vertices $\{x,y\}$ such that $a \le \delta (x,y) \le b$; a gadget, in other words, is just a pair of vertices whose degree difference we can control. 

Once we have found sufficiently many suitable gadgets in $H$, we construct a splitting of $H$ as follows: we use Proposition~\ref{subedge13-approx} to decide, one-by-one for each gadget, which way round to assign the vertices of the gadget to the sides of the splitting. The following lemma makes this idea precise.

\begin{lemma}
\label{subedge13-lemgadget}
Let $H$ be a graph on an even number of vertices and suppose that we can partition $V(H)$ into disjoint collections of pairs $\C{P}_1, \C{P}_2, \dots, \C{P}_k$ such that the pairs in $\C{P}_i$ are $[a_i, b_i]$-gadgets, where $0 \le a_1 \le b_1$ and $0 < a_i \le b_i$ for $2 \le i \le k$. If $b_{i-1} \le a_{i}|\C{P}_{i}|$ for each $2 \le i \le k$, then $V(H)$ can be partitioned into two sets $A, B$ of the same size such that $|\Deg(A)-\Deg(B)| \le b_k$. In particular, if $b_k = 1$, then $H$ is splittable.
\end{lemma}
\begin{proof}
We show by induction on $i$ that it is possible to partition the vertices of the gadgets in $\C{P}_1, \dots, \C{P}_i$ into two sets $A_i$ and $B_i$ of equal size such that $|\Deg(A_i)-\Deg(B_i)| \le b_{i}$. The lemma follows by taking $A=A_k$ and $B=B_k$.

We set $b_0 = 0$ and $A_0 = B_0 = \emptyset$, so the claim is trivially true when $i=0$. So suppose that $i\ge 1$ and that we have constructed $A_{i-1}$ and $B_{i-1}$. Denote the $[a_i, b_i]$-gadgets in $\C{P}_i$ by $(x_j, y_j)$, where $\Deg(x_j) \ge \Deg(y_j)$ for $1 \le j \le |\C{P}_i|$. Using the fact that $b_{i-1} \le a_{i}|\C{P}_{i}|$, it follows from Proposition~\ref{subedge13-approx} that there is a choice of signs $\zeta_j \in \lbrace -1,+1\rbrace$ for $1 \le j \le |\C{P}_i|$ such that 
\[ \big| (\Deg(A_{i-1}) - \Deg(B_{i-1})) + \sum_j \zeta_j \delta(x_j, y_j) \big| \le  b_i.\]
Given $\zeta_j$ as above, we construct $A_i$ and $B_i$ from $A_{i-1}$ and $B_{i-1}$ as follows: for each $1\le j \le |\C{P}_{i}|$, we add $x_j$ to $A_{i-1}$ and $y_j$ to $B_{i-1}$ if $\zeta_j = 1$, and $y_j$ to $A_{i-1}$ and $x_j$ to $B_{i-1}$ if $\zeta_j = -1$. The claim follows.

If $b_k = 1$, notice that we have a partition of $V(H)$ into two sets $A$ and $B$ of equal size such that $|\Deg(A) - \Deg(B)| \le 1$. As $\Deg(A) + \Deg(B)$ is the sum of all the vertex degrees, we conclude that $\Deg(A) = \Deg(B)$ since $\Deg(A) - \Deg(B)$ must be even. 
\end{proof}

Lemma~\ref{subedge13-lemgadget} tells us that a graph is splittable if we can find the right gadgets in the graph. The majority of the work in proving Theorem~\ref{subedge13-mainthm} is in showing that it is possible to find a good collection of gadgets.

\section{Proof of the main result}\label{subedge13-proof}

We now try and make the intuition presented in Section~\ref{subedge13-strat} precise. We shall show that if $\eps > 0$ and $n$ is sufficiently large, it is possible to transform any graph $G$ on $n$ vertices into a splittable graph by deleting at most $2\eps n$ vertices from $G$. Before we begin, we remark that the various constants suppressed by the asymptotic notation throughout the proof are allowed to depend on $\eps$. We shall use $c_1, c_2, \dots$ to represent small constants depending on $\eps$ and $C_1, C_2, \dots$ for large constants depending on $\eps$. All our estimates will hold when $n$ is sufficiently large.

\begin{proof}[Proof of Theorem~\ref{subedge13-mainthm}]

Let $\eps > 0$ be fixed. By deleting an arbitrary vertex of $G$ if necessary, assume that $n=|V(G)|$ is even. Let $\beta = \beta (\eps)$ be a small constant whose value we shall fix at the end of the argument in Case 1.

Call a pair of vertices $\{x,y\}$ a `large' pair if $\delta(x,y) \in [n^{1/3}, \beta n]$. Let $c_1 = \eps / 2 $. We distinguish two cases depending on how many disjoint large pairs we can find in $G$. We first deal with the case when $G$ contains many disjoint large pairs.

\textbf{Case 1: $G$ contains $c_1 n$ disjoint large pairs of vertices.}\label{subedge13-Case1}
In this case, we shall show that $G$ has an induced subgraph $H$ of even order on at least $(1-2\eps)n$ vertices that contains

\begin{enumerate}
\item a collection $\C{S}_H$ of $[1,1]$-gadgets of size $\Omega (n/\log{n})$,
\item a collection $\C{M}_H$ of $[1,n^{2/3}]$-gadgets of size at least $2\beta n$, and
\item a collection $\C{L}_H$ of $[n^{1/9}, 2\beta n]$-gadgets of size $\Omega(n)$
\end{enumerate}
such that the collections $\C{S}_H$, $\C{M}_H$, and $\C{L}_H$ are disjoint. It is straightforward to check that such a graph $H$ is splittable using Lemma~\ref{subedge13-lemgadget}. Indeed, pair up the vertices $V(H) \setminus (L_H\cup M_H \cup S_H)$ arbitrarily; any such pair is a $[0,n]$-gadget, so we have a partition of $V(H)$ into disjoint collections of $[0,n]$-gadgets, $[n^{1/9}, 2\beta n]$-gadgets, $[1,n^{2/3}]$-gadgets and $[1,1]$-gadgets. The sizes of these collections satisfy the conditions of Lemma~\ref{subedge13-lemgadget} if $n$ is sufficiently large and it follows that $H$ is splittable.

We shall now show that $G$ does indeed contain such an induced subgraph $H$. We shall construct $H$ by deleting vertices from $G$ at random.

To avoid notational clutter, in the rest of the argument in Case 1, we shall write \emph{large-gadget} for an $[n^{1/9}, 2\beta n]$-gadget, \emph{medium-gadget} for a $[1,n^{2/3}]$-gadget and \emph{one-gadget} for a $[1,1]$-gadget.

Let $\C{L}$ be a collection of $c_1 n$ large pairs of vertices of $G$. The pairs in $\C{L}$ will be the candidates for the large-gadgets we hope to find in $H$. Our next task is to find a large collection $\C{M}$ of `medium' pairs and a reasonably large collection $\C{S}$ of `small' pairs; the collections $\C{M}$ and $\C{S}$ will provide the candidate pairs for the medium-gadgets and one-gadgets that we would like to find in $H$.

Now, $|V\setminus L| = (1-2c_1)n$; recall that in our notation, $L$ denotes the underlying ground set of $\C{L}$. If we find more than $(1/2 - \eps)n$ disjoint clone pairs $\{x,y\}$ in $G[V\setminus L]$, we are done. Indeed, we can delete all the other ($\le 2\eps n$) vertices not in any of these clone pairs to get a splittable graph: we split this graph by assigning different vertices of each clone pair to different halves of the partition. So we may assume that we can find a set $V' \subset V \setminus L$ of vertices of $G$ such that any two vertices of $V'$ disagree on some vertex of $V\setminus L$ and $|V'| \ge (2\eps - 2c_1)n \ge \eps n$.

Let $C_1 = 4/\eps$ and let $c_2 = \eps / 12$. We now apply Proposition~\ref{subedge13-grouping} to the degrees of the vertices of $V'$; by our choice of $C_1$ and $c_2$, we see that we can find $c_2 n$ disjoint groups of three vertices from $V'$ such that $\delta(x,y) \le C_1$ for any two vertices $x$ and $y$  in the same group. By Proposition~\ref{subedge13-triples}, from each of these triples, we may choose a pair of vertices $\{x,y\}$ such that $\Delta(x,y) \le 2n/3$. Write $\C{P}$ for this collection of $c_2 n$ pairs.

For $0 \le i \le \log{n} - 1$, let $\C{P}_i$ be the collection of those pairs $\{x,y\}$ in $\C{P}$ such that $\Delta (x,y) \in [2^i, 2^{i+1})$. There are two possibilities that we need to consider. It might be that no collection $\C{P}_i$ contains too many pairs; we deal with this case next. The case where one of these collections contains many pairs is easier; we deal with this scenario later with a modification of the argument that follows.

Let $C_2\ge 4$ be a (large) constant depending on $\eps$; we shall fix the value of $C_2$ later in the proof at the end of Case 1A. Also, let $c_3 = c_2 / 3 C_2 \le c_2/ 12$.

\textbf{Case 1A: None of the collections $\C{P}_0, \C{P}_1, \dots, \C{P}_{\log{n}-1}$ contains $c_3n$ pairs.}\label{subedge13-smallbox} It is clear that at least one of the collections $\C{P}_0, \C{P}_1, \dots, \C{P}_{\log{n}-1}$ contains at least $c_2 n / \log{n} $ pairs. Let $k$ be the smallest index such that $|\C{P}_k| \ge c_3 n/ \log{n}$ and let us define our collection of small pairs $\C{S}$ by setting $\C{S} = \C{P}_k$. We now define our collection of medium pairs $\C{M}$ by setting 
\[\C{M} = \C{P}_{k + C_2} \cup \dots \cup \C{P}_{\log{n}-1}.\] 
Since $k$ is minimal and $c_3 \le c_2/12$, we see that $|\C{M}| \ge c_2 n / 2$.

We shall now restrict our attention to the collections $\C{S}$, $\C{M}$ and $\C{L}$; note that they are disjoint. We shall make use of the following facts about these collections.

\begin{enumerate}
\item $\C{S}$ contains $c_3 n / \log{n}$ pairs of vertices $\{x,y\}$ with $\delta(x,y) \le C_1$, $\Delta (x,y) \in[2^k,2^{k + 1})$, and $\Delta (x,y) \le 2n/3$.
\item $\C{M}$ contains $c_2 n/2$ pairs of vertices $\{x,y\}$ such that $\delta(x,y) \le C_1$, and $\Delta (x,y) \ge 2^{k + C_2}$.
\item $\C{L}$ contains $c_1 n$ pairs of vertices $\{x,y\}$ with $\delta(x,y) \in [n^{1/3}, \beta n$].
\item For any pair of vertices $\{x,y\}$ in $\C{S}$ or $\C{M}$, there exists at least once vertex in $V\setminus L$ on which $x$ and $y$ disagree.
\end{enumerate}

We are now in a position to describe how we intend to construct a splittable graph from $G$. We shall delete vertices from $G$ independently with a fixed probability. We shall show that with positive probability, many of the small pairs from $\C{S}$ form one-gadgets in the resulting graph, many of the medium pairs from $\C{M}$ form medium-gadgets, and many of the large pairs from $\C{L}$ form large-gadgets in the resulting graph.

Fix $p = \min \lbrace \eps, 2^{- k} \rbrace$. We now delete vertices from $G$ independently with probability $p$. Let $H$ be the resulting graph. We shall show that with probability $\Omega(1)$, the graph $H$ is splittable and contains at least $(1-2\eps)n$ vertices; this clearly implies the result we are trying to prove.

Note that for a graph to be splittable, it must necessarily contain an even number of vertices. With this in mind, let $\Scr{E}$ be the event that an even number of vertices have been deleted, in other words, $\Scr{E}$ is the event that $|V(H)|$ is even. By Proposition~\ref{subedge13-binparity}, we see that $\P (\Scr{E}) \ge 1/2$. We now analyse what happens to the degree differences of the pairs in $\C{S}$, $\C{M}$ and $\C{L}$ in the graph $H$.

\textsc{One-gadgets.} We first show that many of the pairs in $\C{S}$ form one-gadgets in $H$.

\begin{lemma}
\label{subedge13-onegad}
For any pair $\{x,y\} \in \C{S} $, 
\[ \P ( \{x,y\} \mbox{ is a one-gadget in $H$ }|\, \Scr{E}) \ge f(\eps, C_1) > 0.\]
\end{lemma}

The crucial fact about Lemma~\ref{subedge13-onegad} is that the lower bound on the probability is independent of $C_2$.

\begin{proof}[Proof of Lemma~\ref{subedge13-onegad}]
Let $A = \Gamma (x) \setminus (\Gamma(y)\cup \lbrace y \rbrace) $ and $B =\Gamma (y) \setminus (\Gamma(x) \cup \lbrace x \rbrace)$. Thus, $\delta(x,y) = ||A| - |B||$ and $\Delta (x,y) = |A| +  |B|$. Note that since $x$ and $y$ disagree on at least one vertex of $V \setminus L$, it cannot be the case that both $A$ and $B$ are empty. Suppose without loss of generality that $|A| \ge |B|$ and that in particular, $A \neq \emptyset$.

Let $E_1$ be the event that both $x$ and $y$  are not deleted, $E_2$ the event that no vertices are deleted from $B$, $E_3$ the event that exactly $|\delta(x,y)-1|$ vertices are deleted from $A$,  and $E_4$ the event that the number of vertices deleted from $V \setminus (A \cup B \cup \lbrace x, y \rbrace)$ has the same parity as $|\delta(x,y)-1|$. It is obvious that the family $\{ E_1, E_2, E_3, E_4\}$ is independent since these events correspond to disjoint sets of vertices, and it is easy to check that 
\[ \P ( \{x,y\}\,\mbox{is a 1-gadget in $H$}\,|\,\Scr{E}) \ge \prod_{i=1}^4 \P (E_i). \]

To complete our proof of the claim, we shall bound the factors on the right one by one. Clearly, $\P (E_1) \ge (1-\eps)^2$. 

We trivially have $|A|, |B| \le 2^{k + 1}$. Furthermore $ |A|, |B| \ge 2^{k-1} - C_1/2$, since $0 \le \delta(x,y) \le C_1$. Also, we know that $ \eps 2^{-k} \le p \le 2^{-k}$. To bound $\P(E_2)$, first note that $p|B| \le 2$. Now, $\P(E_2) = \P( \Bi (|B|, p) = 0)$, so by Proposition~\ref{subedge13-binsmall}, $\P(E_2) \ge \exp{(-4)}$.

We now bound $\P( E_3)$.  Clearly, $p|A| \le 2$. If $2^{k} \ge 2C_1$, then $|A| \ge 2^{k-2}$, so $p|A| \ge \eps/4$. If $2^{k} \le 2C_1$, then $p \ge \eps2^{-k} \ge \eps / 2C_1$, so $p|A| \ge \eps/2C_1$ since $|A| \ge 1$. Consequently, 
\[ \min \lbrace \eps/4, \eps / 2C_1 \rbrace \le p |A| \le 2.\] 
Now, $\P(E_3) = \P ( \Bi (|A|, p) = |\delta(x,y)-1|)$. Using the above estimates for $p|A|$ and the fact that $0 \le \delta(x,y) \le C_1$ in Proposition~\ref{subedge13-binsmall}, we see that $\P(E_3) = \Omega_{\eps, C_1} (1)$.

Finally, since $\Delta (x,y) \le 2n/3$, it follows that $|V \setminus (A \cup B)|\ge n/3$, so by Proposition~\ref{subedge13-binparity}, $\P (E_4) \ge 1/6 $ for all sufficiently large $n$. The claim follows.
\end{proof}

From Lemma~\ref{subedge13-onegad} and Proposition~\ref{subedge13-spacesplit} we see that, conditional on $\Scr{E}$, the number of one-gadgets in $H$ from $\C{S}$ is $\Omega ( n / \log{n})$ with probability at least $f(\eps, C_1)/2$; furthermore, and crucially, we note that this lower bound on the probability is independent of $C_2$.

\textsc{Medium-gadgets.}
We next shift our attention to the pairs in $\C{M}$.
\begin{lemma}
\label{subedge13-medgad}
For any pair $\{x,y\} \in \C{M}$, 
\[ \P (1\le \delta(x,y,H) \le n^{2/3}\,|\,x,y \in V(H)) = 1 - o_{C_2 \to \infty}(1) - o(1).
\]
\end{lemma}

\begin{proof}
Let $N_1 = |\Gamma (x) \setminus (\Gamma(y)\cup \lbrace y \rbrace) |$ and let $N_2 =  |\Gamma (y) \setminus (\Gamma(x) \cup \lbrace x \rbrace)|$ and suppose without loss of generality that $N_1 \ge N_2$. Note that $\delta(x,y) = |N_1 - N_2| \le C_1$. Let $X_1$ and $X_2$ be independent random variables with distributions $\Bi (N_1, 1-p)$ and $\Bi (N_2, 1-p)$ respectively. Observe that $\delta(x,y,H)$ has the same distribution as $|X_1 - X_2|$.

We condition on $x,y \in V(H)$. Let $E_1$ be the event that $\delta(x,y,H) = 0$. Clearly, $\P (E_1) = \P (X_1 = X_2)$. Let $E_2$ denote the event that $\delta(x,y,H) \ge n^{2/3}$. It is enough to show that $\P(E_1 \cup E_2) = o_{C_2 \to \infty}(1) + o(1)$.

For any fixed values of $p$ and $N_2$, it is not hard to check that $\P (X_1 = X_2)$ attains its maximum when $N_1 = N_2$; clearly, $\P (X_1 = X_2) = \sum_{i = 0}^{N_2} \P(X_1 = i) \P(X_2 = i)$ and the required conclusion follows from Cauchy--Schwarz inequality. Thus $\P (E_1)$ is bounded above by the probability of two independent random variables with the distribution $\Bi (N_2, 1-p)$, or equivalently $\Bi (N_2, p)$, being equal. Now, $N_2 \ge 2^{k + C_2-1} - C_1/2$ and $p \ge \eps 2^{-k}$. So, $pN_2 \ge \eps 2^{C_2-1} - 2^{-k+1}$ which, since $k \ge 0$, means that $pN_2 \ge \eps 2^{C_2-1} - 2$. As $\eps$ is fixed, we note that $pN_2$ can be made arbitrarily large by choosing $C_2$ large enough. Since $p \le 1/2$, by Proposition~\ref{subedge13-binequal}, we see that $\P (E_1) = o_{C_2 \to \infty}(1)$.

Clearly, $\P(E_2) = \P(|X_1 - X_2| \ge n^{2/3} )$. Applying Proposition~\ref{subedge13-binfar} to $X_1$ and $X_2$, we conclude that $\P (E_2) = O (\exp ( -n^{1/3}/5))$.
\end{proof}

Let $\C{M}'$ be the collection of those pairs $\{x,y\}\in \C{M}$ such that both $x$ and $y$ survive in $H$. Since the family of events $\lbrace x, y \in V(H) \rbrace$ is a family of mutually independent events for different pairs $\{x,y\}\in \C{M}$ and since $\P (x, y \in V(H)) \ge (1-\eps)^2$, it follows from Proposition~\ref{subedge13-bernstein} that $\P (|\C{M}'| < (1-\eps)^2 |\C{M}|/2 ) = \exp (-\Omega (n))$.

Consequently, from Lemma~\ref{subedge13-medgad}, it follows that for any pair $\{x,y\} \in \C{M}$,
\begin{multline*}
\P \left(1\le \delta(x,y,H) \le n^{2/3}\,\Big| \, \{x,y\} \in \C{M}', \, |\C{M}'|>\frac{(1-\eps)^2 |\C{M}|}{2} \right) \\
= 1 - o_{C_2 \to \infty}(1) - o(1).
\end{multline*} 
Thus by Proposition~\ref{subedge13-spacesplit}, conditional on $|\C{M}'|>(1-\eps)^2 |\C{M}|/2$, the number of medium-gadgets from $\C{M}'$ in $H$ is at least $|\C{M}'|/3$ with probability $1 - o_{C_2 \to \infty}(1) - o(1)$. Thus, the number of medium-gadgets in $H$ is at least $(1-\eps)^2 |\C{M}| /6$ with probability $(1 - \exp (-\Omega(n)))(1 -  o_{C_2 \to \infty}(1) - o(1))$, which is still $1 -  o_{C_2 \to \infty}(1) - o(1)$.

Thus, conditional on the event $\Scr{E}$, the number of medium-gadgets in $H$ from $\C{M}$ is $\Omega(n)$ with probability $1 -  o_{C_2 \to \infty}(1) - o(1)$.

\textsc{Large-gadgets.}
We finally consider the pairs of vertices in $\C{L}$. Recall that every pair $\{x,y\} \in \C{L}$ is such that $\delta(x,y) \in [n^{1/3}, \beta n]$, where $\beta$ is a (small) constant whose value we have yet to fix. (Indeed, the value of $\beta$ has so far played no role in our calculations.)

\begin{lemma}
\label{subedge13-largegad}
For any pair $\{x,y\} \in \C{L}$, 
\[\P (n^{1/9} \le \delta(x,y,H) \le 2\beta n\,|\,x,y \in V(H)) = 1 - o(1). \]
\end{lemma}
\begin{proof}
We condition on $x,y \in V(H)$. Let $E_1$ be the event that $\delta(x, y, H) < n^{1/9}$. Since $\delta(x,y) \ge n^{1/3}$, it follows immediately from Proposition~\ref{subedge13-binclose} that $\P (E_1) = o(1)$.

Let $E_2$ be the event that $\delta(x, y, H) > 2\beta n$. Let $A = \Gamma (x) \setminus (\Gamma(y)\cup \lbrace y \rbrace) $ and $B =\Gamma (y) \setminus (\Gamma(x) \cup \lbrace x \rbrace)$, and let $X_1$ and $X_2$ be random variables that denote the the number of vertices from $A$ and $B$ respectively which survive in $H$. Clearly, the distributions of $X_1$ and $X_2$ are $\Bi (|A|, 1-p)$ and $\Bi (|B|, 1-p)$ respectively.

If $E_2$ were to occur, i.e., it were the case that $|X_1 - X_2| > 2\beta n $, then this would imply that either $|X_1 - (1-p)|A|| \ge \beta n /2$ or $|X_2 - (1-p)|B|| \ge \beta n /2$, since $(1-p)||A| - |B|| \le \delta(x,y) \le \beta n$. It follows that $\P (E_2) = o(1)$ since the probability of either of the above two possibilities is $\exp {(-\Omega(n))}$ by Proposition~\ref{subedge13-bernstein}.
\end{proof}

Arguing as in the case of medium-gadgets, we see from Lemma~\ref{subedge13-largegad} that conditional on the event $\Scr{E}$, the number of large-gadgets in $H$ from $\C{L}$ is $\Omega (n)$ with probability $1 - o(1)$.

\textsc{Constructing a splitting.} We now have a reasonably clear picture of what the degree differences in $H$ of the pairs of vertices in $\C{S}$, $\C{M}$ and $\C{L}$ look like. In summary, conditional on $\Scr{E}$, we have demonstrated that in $H$, we can find

\begin{enumerate}
\item a collection $\C{S}_H$ of $\Omega (n/\log{n})$ one-gadgets with probability $f(\eps, C_1)/2$,
\item a collection $\C{M}_H$ of $\Omega(n)$ medium-gadgets with probability $1-o_{C_2 \to \infty}(1) -o(1)$, and
\item a collection $\C{L}_H$ of $\Omega(n)$ large-gadgets with probability $1-o(1)$
\end{enumerate}
such that the collections $\C{S}_H, \C{M}_H$ and $\C{L}_H$ are disjoint.

Thus by choosing $C_2$ to be a sufficiently large constant depending on $\eps$, by the union bound, we find all of the above with probability $\Omega (1)$ conditional on $\Scr{E}$, provided $n$ is sufficiently large. Also, the expected number of vertices deleted is at most $\eps n$, so by Proposition~\ref{subedge13-bernstein}, the probability that we have deleted more than $2\eps n $ vertices is $\exp {(-\Omega(n))}$.

Consequently, we see that $H$, with probability $\Omega (1)$, has the aforementioned collections of gadgets, and furthermore, also has an even number of vertices and at least $(1-2\eps) n$ vertices. We are done if we can guarantee that $2\beta n \le |\C{M}_H|$; this is possible if we choose $\beta=\beta(\eps)$ to be a suitably small constant because $|\C{M}_H| = \Omega(n)$.

We now consider the case where one of the sets $\C{P}_i$ contains many pairs.

\textbf{Case 1B: One of the sets $\C{P}_0, \C{P}_1, \dots, \C{P}_{\log{n}-1}$ contains $c_3n$ pairs.}\label{subedge13-largebox} This case is easier to deal with than the previous one. We shall argue exactly as before; however we shall have no need of medium-gadgets and it will suffice to consider one-gadgets and large-gadgets alone.

Let $k$ be any index such that $|\C{P}_k|\ge c_3 n$ (while we chose $k$ to be minimal previously, any index $k$ such that $|\C{P}_k|\ge c_3 n$ will do in this case). As before, we set $p = \min \lbrace \eps, 2^{- k} \rbrace$ and $\C{S}=\C{P}_k$. We now delete vertices from $G$ independently with probability $p$. Let $H$ be the resulting graph; as before, we condition on deleting an even number of vertices. We claim that $H$ is splittable with probability $\Omega (1)$.

It is not hard to check that Lemma~\ref{subedge13-onegad} and Lemma~\ref{subedge13-largegad} hold in this case as well. We conclude that we can delete an even number of vertices from $G$ to obtain a graph $H$ with $|V(H)| \ge (1-2\eps) n$ in such a way that in $H$, we can find
\begin{enumerate}
\item a collection $\C{S}_H$ of $\Omega (n)$ one-gadgets, and
\item a collection $\C{L}_H$ of $\Omega(n)$ large-gadgets
\end{enumerate}
such that $\C{S}_H$ and $\C{L}_H$ are disjoint. As before, it follows from Lemma~\ref{subedge13-lemgadget} that $H$ is splittable when $n$ is sufficiently large provided $2\beta n \le |\C{S}_H|$; this is possible if we choose $\beta=\beta(\eps)$ to be a suitably small constant because $|\C{S}_H| = \Omega(n)$.

Thus, for all sufficiently small $\beta$ (so as to satisfy the conditions from both Case 1A and 1B), we see that we are done if $G$ contains many disjoint large pairs. Note that we have now fixed the value of $\beta$. We now deal with the case $G$ does not contain many disjoint large pairs.

\textbf{Case 2: $G$ does not contain $c_1 n$ disjoint large pairs.}\label{subedge13-Case2}
In this case, we shall show that $G$ has an induced subgraph $H$ of even order on at least $(1-2\eps)n$ vertices such that $V(H)$ may be partitioned into

\begin{enumerate}
\item a collection $\C{S}_H$ of $[1,1]$-gadgets of size $\Omega (n/\log{n})$, and
\item a collection $\C{M}_H$ of $[0,n^{2/3}]$-gadgets.
\end{enumerate}

In the rest of the argument in Case 2, we shall, as before, call $[1,1]$-gadgets \emph{one-gadgets} and we call $[0,n^{2/3}]$-gadgets (as opposed to $[1,n^{2/3}]$-gadgets as we did earlier) \emph{medium-gadgets}.
 
It is easily seen from Lemma~\ref{subedge13-lemgadget} that $H$ is splittable if $n$ is sufficiently large. We construct our splitting by starting with the pairs in $\C{M}_H$ - we can use these pairs to construct a partition such that sums of the degrees of the vertices of the two halves of the partition differ by at most $n^{2/3}$. We then use the the pairs in $\C{S}_H$ to reduce the difference to at most one; we are done by parity considerations.

We now show how to find such a subgraph $H$. We start by describing how to find pairs of vertices which will be the candidates for the medium-gadgets we hope to find in $H$.

Let $\C{L}$ be a maximal collection of large pairs in $G$. Note that since $\C{L}$ is maximal, we have either $\delta(x,y) < n^{1/3}$ or $\delta(x,y) > \beta n$ for any two vertices $x,y \in V \setminus L$. As $\beta n > 2 n ^ {1/3}$ for all sufficiently large $n$, there is a partition $V \setminus L = K_1 \cup K_2 \cup \dots \cup K_m$ into `clumps' $K_i$ with $m \le 1/\beta$ in such a way that $\delta(x,y) < n^{1/3}$ for any $ x,y \in K_i$ and $\delta(x,y) > \beta n $ if $x\in K_i$ and $y \in K_j$ with $i \neq j$.

We ignore the way in which vertices are originally paired in $\C{L}$ and focus on the ground set $L$. By Proposition~\ref{subedge13-grouping}, we can find from $L$, at least $|L|/2 - n^{1/2}$ disjoint pairs $\{x,y\}$ such that $\delta(x,y) \le n^{1/2}$; call this collection of pairs $\C{Q}$.

Let $F$ be the graph obtained from $G$ as follows. Delete every vertex of $L \setminus Q$. Delete one vertex from every clump $K$ which contains an odd number of vertices. Having done this, delete a clump $K$ (i.e., delete all the vertices of $K$) if $|K| \le n^{1/2}$.

Note that the vertex set of $F$ consists of the surviving clumps, each of which has even size and cardinality at least $n^{1/2}$, and the (possibly empty) set of pairs $\C{Q}$. Since we had at most $1/\beta$ clumps initially, we have deleted $O (n^{1/2})$ vertices in total from $G$ to obtain $F$. Hence, for any two vertices $x,y\in V(F)$, $|\delta(x,y,F) - \delta(x,y,G)| = O(n^{1/2})$. Hence, if either $x$ and $y$ both belong to the same (surviving) clump or if the pair $\{x,y\}$ is in $\C{Q}$, then $\delta(x,y,F) = O (n^{1/2})$. Let us say that two vertices $x,y \in V(F)$ are \emph{proximate} if either both $x$ and $y$  belong to the same clump in $F$ or if $\{x,y\}\in \C{Q}$; these proximate pairs of vertices will be our candidates for medium-gadgets in $H$.

We now show how to find pairs of vertices which will be the candidates for the one-gadgets we hope to find in $H$. We shall henceforth work with $F$ as opposed to $G$. We shall write $V$ for $V(F)$ and  all degrees and degree differences, unless specified otherwise, will be with respect to $F$.

Since $|\C{L}|\le c_1 n = \eps n /2$ and since we have only deleted $O (n^{1/2})$ vertices so far, note that $|V \setminus Q| \ge (1 - 3\eps / 2)n$ for $n$ sufficiently large.

If we find at least $(1/2 - \eps)n$ disjoint clone pairs $\{x,y\}$ in $F[V\setminus Q]$, we are done. So we may assume that we can find a set $V' \subset V \setminus Q$ of vertices of $F$ with $|V'| \ge (2\eps - 3\eps/2)n = \eps n / 2$ such that any two vertices of $V'$ disagree on some vertex in $V\setminus Q$.

We claim that if $C_3$ is sufficiently large (as a function of $\beta$), then we can find from any subset of $C_3$ vertices of $V'$, two vertices $x$ and $y$  such that for each clump $K$, the number of vertices of $K$ on which $x$ and $y$  disagree is at most $2|K|/3$. To see this, suppose that we have found $C_3$ vertices such that any two of them $x$ and $y$  disagree on more than two thirds of some clump $K_{x,y}$. Applying Ramsey's theorem (with $1/\beta$ colours) to the complete graph on these $C_3$ vertices where the edge between $x$ and $y$  is labelled by the clump $K_{x,y}$, we see that we can find a monochromatic triangle provided $C_3$ is large enough. But by Proposition~\ref{subedge13-triples}, out of any three vertices, at least two disagree on at most two thirds of the vertices of $K$. We have a contradiction.

Choose $C_3$ as described above and set $C_4 = 4C_3 / \eps$ and $c_4 = \beta / 2C_4$. By Proposition~\ref{subedge13-grouping}, we can find from $V'$, at least $n/C_4$ disjoint groups of size $C_3$ such that that $\delta(x,y) \le C_4$ for any two vertices $x$ and $y$ in the same group. From each of these $n/C_4$ groups of size $C_3$, choose a pair of vertices $\{x,y\}$ such that $x$ and $y$  disagree on at most two thirds of every clump. Choose a clump $K^*$ such that at least a $\beta$ fraction of these pairs $\{x,y\}$ are such that $x$ and $y$  disagree on at least one vertex in $K^*$; this is possible because any two vertices of $V'$ disagree on $V(F) \setminus Q$ and consequently, on at least one clump and furthermore, there are at most $1/\beta$ clumps. Let $\C{P}$ be this collection of pairs which all disagree on at least one vertex in $K^*$; clearly, $|\C{P}| \ge  \beta n / C_4 = 2c_4n$.

We shall proceed as in Case 1 by pigeonholing the pairs in $\C{P}$ into different boxes based on the size of their difference neighbourhoods, but with one important difference. Note that while \emph{any two} vertices in the same clump have a small ($O(n^{1/2})$) degree difference, we can only guarantee that two vertices of $Q$ have small ($O(n^{1/2})$) degree difference if the pair belongs to $\C{Q}$. Consequently, when we later delete vertices at random, we shall either delete both vertices of a pair in $\C{Q}$ or retain both; hence, we shall treat a pair of vertices in $\C{Q}$ as a single vertex when it comes to pigeonholing the pairs in $\C{P}$. This is made precise below.

Let $F_\C{Q}$ be the multigraph without loops obtained from $F$ by contracting every pair $\{x,y\}$ in $\C{Q}$ (we ignore the loops that might arise). Note that there are at most two parallel edges between any two vertices of $F_\C{Q}$ unless both vertices correspond to contracted pairs from $\C{Q}$, in which case there are at most four parallel edges between them. In $F_\C{Q}$, we say that two vertices $x$ and $y$  disagree on a vertex $v \neq x,y$ if the number of edges between $v$ and $x$ is not equal to the number of edges between $v$ and $y$. For $0 \le i \le \log{n}-1$, let $\C{P}_i$ be the collection of those pairs $\{x,y\}$ in $\C{P}$ such that $\Delta (x,y,F_\C{Q}) \in [2^i, 2^{i+1})$, where $\Delta(x,y,F_\C{Q})$ is the number of vertices of $F_\C{Q}$ on which $x$ and $y$  disagree. 

As before, let $k$ be any index such that $|\C{P}_k|\ge 2c_4 n / \log{n}$; take $\C{S} = \C{P}_k$ and set $p = \min \lbrace \eps, 2^{- k} \rbrace$.

In summary, $\C{S}$ consists of pairs $\{x,y\}$ such that 
\begin{enumerate}
\item $x$ and $y$  disagree on at most two thirds of every clump,
\item $x$ and $y$  disagree on at least one vertex of $K^*$, 
\item $\delta(x,y)\le C_4$, and 
\item $\Delta(x,y,F_\C{Q}) \in [2^k, 2^{k +1})$. 
\end{enumerate}
Furthermore, since $\delta(x,y)\le C_4 =o(n^{1/2})$ for any $\{x,y\}\in \C{S}$, both members of any pair in $\C{S}$ must belong to the same clump.

Consider the partition $\C{S} = \C{S}_o \cup \C{S}_e$ where $\C{S}_o$ is the set of those pairs $\{x,y\}\in \C{S}$ such that $\delta(x,y)$ is odd. Recall that $|\C{S}| \ge 2c_4 n / \log{n}$, so one of $\C{S}_o$ or $\C{S}_e$ contains more than $c_4n/\log{n}$ pairs. At this point, we need slightly different arguments depending on whether we have more pairs with odd degree difference or even degree difference in $\C{S}$.

\textbf{Case 2A: $\C{S}$ contains many odd pairs.}
We first consider the case where $|\C{S}_o| \ge c_4 n / \log{n}$. We shall delete vertices from $F$ as follows. We pick vertices of $F_\C{Q}$ independently with probability $p$. For every vertex of $F_\C{Q}$ that we pick, we delete (as appropriate) either the corresponding vertex or both vertices of the corresponding pair of vertices from $\C{Q}$ in $F_\C{Q}$. Let $H$ be the resulting graph. Our aim is to show that $H$ is splittable with probability $\Omega (1)$.

Earlier, we conditioned on deleting an even number of vertices from $G$. In this case, we need a little more. Let $\Scr{E}^*$ be the event that an even number of vertices were deleted from each clump. By Proposition~\ref{subedge13-binparity}, we see that $\P (\Scr{E}^*) \ge (1/2)^{1/ \beta}$. Note that a consequence of $\Scr{E}^{*}$ is that $|V(H)|$ is even.

\textsc{One-gadgets.}
First, we shall show that many of the pairs in $\C{S}_o$ become one-gadgets in $H$.

\begin{lemma}
\label{subedge13-onegad2}
For any pair $\{x,y\} \in \C{S}_o $, 
\[ \P ( \{x,y\} \mbox{ is a one-gadget in $H$ } | \,\Scr{E}^{*}) =  \Omega (1).\]
\end{lemma}
\begin{proof}
In $F_\C{Q}$, let $A$ be the set of those vertices $v\neq x,y$ such that number of edges between $v$ and $x$ is more than the number of edges between $v$ and $y$ and let $B$ be defined analogously by interchanging $x$ and $y$. Let $A = A_1 \cup A_2$ where $A_1$ and $A_2$ are respectively those vertices $v$ in $A$ such that the number of edges between $v$ and $x$ is one, respectively two, more than the number of edges from $v$ to $y$; define $B_1$ and $B_2$ analogously.

The proof follows that of Lemma~\ref{subedge13-onegad}. Clearly,
\[ 2^k \le |A_1| + |A_2| + |B_1| + |B_2| < 2^{k +1}.\]
Furthermore, $\delta(x,y) = ||A_1| + 2|A_2| - |B_1| - 2|B_2||$, so 
\[ -C_4 \le |A_1| + 2|A_2| - |B_1| - 2|B_2| \le C_4.\] 
Using the above two inequalities, it is not hard to check that 
\[\max \lbrace |A_1|, |A_2| \rbrace, \max \lbrace |B_1|, |B_2| \rbrace \ge 2^{k - 3} - C_4/4.\]

Since $\delta(x,y)$ is odd, suppose without loss of generality that $\Deg(x)>\Deg(y)$. Let $E_1$ be the event that both $x$ and $y$  are not picked to be deleted, $E_2$ the event that no vertices are picked from $B$, $E_3$ the event that $X_1 + 2X_2 = \delta(x,y)-1$ where $X_1$ and $X_2$ are the number of vertices picked from $A_1$ and $A_2$ respectively, and $E_4$ the event that the number of vertices picked from $K \setminus (A \cup B \cup \lbrace x, y \rbrace)$ has the same parity as the number of vertices picked from $K \cap (A \cup B \cup \lbrace x, y \rbrace)$ for every clump $K$. The collection of events $\{E_1,E_2,E_3\}$ is clearly independent, and it is easy to check that 
\[ \P ( \{x,y\} \mbox{ is a one-gadget } | \, \Scr{E}^{*}) \ge  \P (E_1) \P (E_2) \P (E_3) \P (E_4|E_1,E_2,E_3).\]

Clearly, $\P (E_1) \ge (1-\eps)^2$. As in Lemma~\ref{subedge13-onegad}, note that $p|B| \le 2$, so by Proposition~\ref{subedge13-binsmall}, $\P (E_2) \ge \exp {(-4)}$.

We now bound $\P (E_3)$. First suppose that $2^{k - 3} - C_4/4 > C_4$. Recall that $\delta(x,y)$ is odd. If $|A_2| \ge |A_1|$, we consider the event that $(\delta(x,y)-1)/2$ vertices are picked from $A_2$ and no vertices are picked from $A_1$ in $F_\C{Q}$; as in Lemma~\ref{subedge13-onegad}, we see that $p|A_2| = \Theta(1)$, so this event occurs with probability $\Omega (1)$. Hence, $E_3$ occurs with probability $\Omega (1)$. If $|A_1| > |A_2|$, we consider the event that $\delta(x,y)-1$ vertices are picked from $|A_1|$ and no vertices are picked from $A_2$ and note that this event occurs with probability $\Omega (1)$; hence, $E_3$ occurs with probability $\Omega (1)$.

If, on the other hand, $2^{k - 3} - C_4/4 \le C_4$, then clearly $k = \Theta (1)$l hence, $p, |A_1|, |A_2|$ are all $\Theta (1)$. In this case, we consider the event that $t = \min \lbrace (\delta(x,y)-1)/2, |A_2| \rbrace$ vertices are picked from $A_2$ and $\delta(x,y)-1-2t$ vertices are picked from $A_1$. Now,  $|A_1| + 2|A_2| \ge \delta(x,y)$ since we assumed that $\Deg(x) > \Deg(y)$, so $|A_1| \ge \delta(x,y)-1-2t$. Also, as noted above, $p, |A_1|, |A_2|$ are all $\Theta (1)$. So this event occurs with probability $\Omega (1)$. Hence, the event $E_3$ occurs with probability $\Omega (1)$.

Since $x$ and $y$  disagree on at most two thirds of every clump and since every clump has size at least $n^{1/2}$, it follows from Proposition~\ref{subedge13-binparity} that $\P (E_4|E_1, E_2, E_3) \ge (1/6)^{1/\beta}$ for all sufficiently large $n$.
\end{proof}

Let $\C{S}_H$ be the set of pairs from $\C{S}_o$ that form one-gadgets in $H$. From Lemma~\ref{subedge13-onegad2} and Proposition~\ref{subedge13-spacesplit}, we see that conditional on $\Scr{E}^{*}$, $|\C{S}_H| \ge \E [|\C{S}_H|]/2 = \Omega ( n / \log{n})$ with probability $\Omega (1)$.

\textsc{Medium-gadgets.}
 We now show that the degree difference of any pair of vertices which are proximate in $F$ cannot become too large in $H$.

\begin{lemma}
\label{subedge13-med2}
Conditional on $\Scr{E}^{*}$ and $|\C{S}_H|\ge \E [|\C{S}_H|]/2$, the probability that there exist $x,y\in V(H)$ which are proximate in $F$ and satisfy $\delta(x,y,H) > n^{2/3}$ is $o(1)$.
\end{lemma}
\begin{proof}
Recall that for any two vertices $x$ and $y$  which are proximate in $F$, $\delta(x,y) = O(n^{1/2})$. For such a pair of  vertices $x$ and $y$, note by Proposition~\ref{subedge13-binfar} that 
\[\P (\delta(x,y,H) > n^{2/3}\,|\,x,y \in V(H)) = O(\exp (  -n^{1/3}/5)).\] 
Consequently, since we have conditioned on an event with probability $\Omega (1)$, the probability that there exist some vertices $x,y \in V(H)$ such that $x$ and $y$ are proximate and $\delta(x,y,H)> n^{2/3}$  is $O(n^2 \exp ( -n^{1/3}/5)) = o(1)$.
\end{proof}

\textsc{Constructing a splitting.} We now describe how to construct a splitting of $H$. Let $\C{Q}_H$ be the set of pairs from $\C{Q}$ that survive in $H$. For a clump $K$ in $F$, let $K_H$ denote the set $(K\setminus S_H)\cap V(H)$. Clearly, $V(H)$ is the disjoint union of $S_H$, $Q_H$ and the clumps $K_H$. Note that conditional on $\Scr{E}^{*}$, the size of $K_H$ is even for every clump $K$ since both members of any pair in $\C{S}_H$ must necessarily belong to the same clump. Since each $K_H$ has even cardinality, we may group the vertices of each $K_H$ into pairs. Pair up the vertices in each $K_H$ arbitrarily; let $\C{M}_H$ be the collection consisting of these pairs and the pairs in $\C{Q}_H$. Clearly, every pair of vertices in $\C{M}_H$ are proximate in $F$ and by Lemma~\ref{subedge13-med2}, the probability that there exists some pair $\{x,y\} \in \C{M}_H$ with $\delta(x,y,H) > n^{2/3}$ is $o(1)$.

The expected number of vertices deleted from $F$ is at most $\eps n$ and the number of vertices deleted from $G$ to obtain $F$ is $O(n^{1/2})$. Hence, by Proposition~\ref{subedge13-bernstein}, the probability that we have deleted more than $2\eps n $ vertices from $G$ is $\exp {(-\Omega(n))}$.

We conclude that there exists an induced subgraph $H$ of $G$ such that $|V(H)|\ge (1-2\eps)n$, and with the further property that $V(H)$ may be partitioned into

\begin{enumerate}
\item a collection $\C{S}_H$ of one-gadgets of size $\Omega (n/\log{n})$, and
\item a collection $\C{M}_H$ of medium-gadgets.
\end{enumerate}

It follows from Lemma~\ref{subedge13-lemgadget} that $H$ is splittable and we are done.

\textbf{Case 2B: $\C{S}$ contains many even pairs.} \label{subedge13-even}
Now we consider the case where $|\C{S}_e| \ge c_4 n / \log{n}$. 

Note that since we intend to delete either both vertices of a pair in $\C{Q}$ or neither, it might be the case that it is impossible to make the parity of the degree difference of a pair in $\C{S}_e$ odd in $H$. Consequently, in this case, we will need to work with $[2,2]$-gadgets, or \emph{two-gadgets} for short, in addition to one-gadgets. With the exception of this slight change of tack to account for parity considerations, the argument is quite similar to the one in the previous case, and we only sketch it.

Let $c_5$ be a (small) constant depending on $\eps$; the value of $c_5$ will be chosen later, following the statement of Lemma~\ref{subedge13-twogad}.

Recall that every pair of vertices in $\C{S}_e$ disagree on some vertex in the clump $K^*$. Suppose there exists a vertex $v\in K^*$ such that $c_5 n/ \log{n}$ pairs from $\C{S}_e$ all disagree on $v$. In this case, we may complete the proof as follows. Let $\C{S}_v\subset \C{S}_e$ be the collection of pairs in $\C{S}_e$ that disagree on $v$. We shall delete vertices from $F$ as follows. We first delete $v$ and then delete one other vertex uniformly at random from $K^*$. Following this, we proceed as before by picking vertices of $F_\C{Q}$ independently with probability $p$ and then deleting the corresponding vertices or pairs of vertices from $\C{Q}$ in $F$. Let $H$ be the resulting graph. Note that when we delete $v$, the degree difference of every pair in $\C{S}_v$ changes parity and becomes odd. When we then delete another vertex uniformly at random from $K^*$, the parity of the degree difference of a pair in $\C{S}_v$ is unaltered with probability at least $1/3$ since every pair in $\C{S}$ disagree on at most two thirds of any clump. Arguing as in Lemma~\ref{subedge13-onegad2}, for any pair in $\C{S}_v$, we see that the probability that this pair forms a one-gadget in $H$, conditional on deleting an even number of vertices from every clump, is $\Omega(1)$ (albeit with a smaller constant than in Case 2A). Since $|\C{S}_v| \ge c_5 n/ \log{n}$, we can conclude the proof exactly as in the case where $\C{S}$ contains many odd pairs.

Hence, we may assume that for every vertex $v\in K^*$, the number of pairs in $\C{S}_e$ that disagree on $v$ is at most $c_5 n/ \log{n}$. We delete vertices from $F$ as before by picking vertices of $F_\C{Q}$ independently with probability $p$ and then deleting the corresponding vertices or pairs of vertices from $\C{Q}$ in $G$. Let $H$ be the resulting graph.

As before, let $\Scr{E}^{*}$ be the event that an even number of vertices were deleted from each clump. The proof of Lemma~\ref{subedge13-onegad2}, with minor modifications for the change in parity, yields a proof of the following lemma.

\begin{lemma}
\label{subedge13-twogad}
For any $\{x,y\} \in \C{S}_e $, $\P ( \{x,y\}\, \mbox{is a two-gadget in $H$}\,|\, \Scr{E}^{*}) =  \Omega (1)$. \qed
\end{lemma}

Let $\C{S}_{H}$ be the collection of pairs from $\C{S}_e$ that form two-gadgets in $H$. From Lemma~\ref{subedge13-twogad} and Proposition~\ref{subedge13-spacesplit}, we see that there exists a small positive constant $c_6$ such that, conditional on $\Scr{E}^{*}$, $|\C{S}_{H}| \ge c_6 n/ \log{n}$ with probability $\Omega (1)$. Let us now fix $c_5 = c_6/4$.

\textsc{Constructing a splitting.}
As before, let $\C{Q}_{H}$ be the collection of pairs from $\C{Q}$ that survive in $H$ and for each clump $K$ in $F$, let us write $K_{H}$ for the set $(K\setminus S_{H})\cap V(H)$.

We have shown that with probability $\Omega(1)$, the graph $H$ is such that
\begin{enumerate}
\item $|K_{H}|$ is even for every clump $K$, and
\item $|\C{S}_{H}| \ge c_6n/ \log{n}$.
\end{enumerate}

Consider any pair $\{x,y\} \in \C{S}_{H}$ and note that in $F$, $x$ and $y$ disagree on at most two thirds of any clump; in particular, $x$ and $y$ agree on at least a third of $K^*$. Consequently, the probability that $x$ and $y$ disagree on every vertex of $K^*_{H}$ is $\exp {(-\Omega(n^{1/2}))}$. Hence, with probability $1-o(1)$, for every $\{x,y\}\in \C{S}_{H}$, there exists some vertex in $K^*_{H}$ on which $x$ and $y$ agree. 

Next, it follows from Lemma~\ref{subedge13-med2} that with probability $1-o(1)$, any two vertices $x,y \in V(H)$ which are proximate satisfy $\delta(x,y,H) \le n^{2/3}$. Finally, the probability that we have deleted more that $2\eps n - 2$ vertices of total from $G$ is, by Proposition~\ref{subedge13-bernstein}, $\exp {(-\Omega(n))}$. It follows that with probability $\Omega(1)$, the graph $H$, in addition to possessing the aforementioned properties, also has the following properties.

\begin{enumerate}
\setcounter{enumi}{2}
\item For every $\{x,y\}\in \C{S}_{H}$, there exists some vertex in $K^*_{H}$ on which $x$ and $y$ agree.
\item For any $x,y\in V(H)$ such that $x$ and $y$ are proximate in $F$, $\delta(x,y,H) \le n^{2/3}$.
\item $ |V(H)| \ge (1-2\eps)n+2$.
\end{enumerate}

With a view of making the graph $H$ splittable, we alter $H$ as follows. Fix a pair $(x^*,y^*)\in \C{S}_{H}$ and a vertex $v\in K^*$ on which $x^*$ and $y^*$ disagree. We know that there is a vertex $u \in K^*_{H}$ on which $x^*$ and $y^*$ agree. Delete $u$ from $H$. If $v \in V(H)$, delete $v$ from $H$ and if $v\notin V(H)$, add $v$ back. After these alterations, note that $H$ still has an even number of vertices. Note also that now, $|V(H)| \ge (1-2\eps)n$ and $\delta(x^*,y^*,H)\in \lbrace 1,3\rbrace$. 

Before we altered $H$, at most $c_5 n/ \log {n}$ pairs in $\C{S}_{H}$ disagreed on any vertex in $K^*$; the alterations above change the degree differences of at most $2c_5 n /\log{n}= c_6 n /2 \log{n}$ pairs in $\C{S}_{H}$. Hence, $H$ contains a collection $\C{S}_H$ of least $c_6 n /2 \log{n} - 1$ pairs of vertices $\{x,y\}$ such that $\delta(x,y,H)=2$ and a pair $(x^*,y^*)$ such that $\delta(x^*,y^*,H)\in \lbrace 1,3\rbrace$. Furthermore, all the vertices of $V(H) \setminus (S_H \cup \lbrace x^*, y^* \rbrace)$ may be grouped into pairs $\{x,y\}$ such that $\delta(x,y,H) \le n^{2/3} + 2$; let $\C{M}_H$ denote this collection of pairs.

It is now easy to check that $H$ is splittable using the argument used to prove Lemma~\ref{subedge13-lemgadget}. Indeed, we can use pairs in $\C{M}_H$ to construct a partition such that sums of the degrees of the vertices of the two halves of the partition differ by at most $n^{2/3} + 2$. For $n$ sufficiently large, we can then reduce the difference to at most two by using all but one of the pairs in $\C{S}_H$. Finally, using the one remaining pair in $\C{S}_H$ and the pair $(x^*,y^*)$, we can reduce the difference to at most one; we are done constructing a splitting of $H$ by parity considerations. This completes the proof of Theorem~\ref{subedge13-mainthm}.
\end{proof}

\section{Conclusion}\label{subedge13-conc}
We have shown that $f(n) \ge n/2 - o(n)$. In fact, it should be possible to read out a bound of $f(n) \ge n/2 - n/(\log{\log{n}})^c$ from our proof for some absolute constant $c>0$; we chose not to include a proof of this fact to keep the presentation simple, and because we do not believe that such a bound is close to the truth. While we have managed to pin down $f$ up to its first order term, there is still a large gap between the upper and lower bounds for $n/2 - f(n)$.

\begin{problem}
What is the asymptotic behaviour of $n/2 - f(n)$?
\end{problem}

We know that $n/2-f(n) = \Omega (\log{\log{n}})$ and $n/2-f(n) = o(n)$; we suspect that the truth lies closer to the lower bound and that in particular, $n/2 - f(n) = o(n^\eps)$ for every $\eps>0$. Indeed, it is not inconceivable that $n/2 - f(n) = \Theta (\log{n})$.

It is natural to generalise the problem to the case where we have more than one type of edge, or ask for more than two disjoint subgraphs. For any $r,l\in \N$, given an edge colouring $\Delta$ of the complete graph on $n$ vertices with $r$ colours, let $g(\Delta)$ be the largest integer $k$ for which we can find $l$ disjoint subsets $V_1,V_2,\dots, V_l$ of $[n]$, each of cardinality $k$, such that for each $1\le i \le r$, the number of edges induced by $V_j$ of colour $i$ is the same for every $1 \le j \le l$. Let $g(n,r,l)$ be the minimum value of $g(\Delta)$ taken over all edge colourings of the complete graph on $n$ vertices. In particular, note that $g(n,2,2) = f(n)$. We conjecture that $g(n,r,2) = n/2 - o(n)$ and more generally, ask the following question.

\begin{problem}
For $r,l\in \N$, what is the asymptotic behaviour of $g(n,r,l)$?
\end{problem}

Finally, we mention a question about digraphs that we find particularly appealing. Given a digraph $D$ on $n$ vertices, let $h(D)$ denote the largest integer $k$ for which there exist disjoint subsets $A,B\subset V$ such that $|A| = |B|=k$ and the number of directed edges from $A$ to $B$ is equal to the number of directed edges from $B$ to $A$. Let $h(n)$ be the minimum value of $h(D)$ taken over all digraphs on $n$ vertices.

\begin{problem}
Determine $h(n)$.
\end{problem}

\section*{Acknowledgements} The first author was partially supported by NSF grant DMS-1301614 and EU MULTIPLEX grant 317532.
 
The majority of the research in this paper was carried out while the authors were visitors at Microsoft Research, Redmond. We are grateful to Yuval Peres and the other members of the Theory Group at Microsoft Research for their hospitality.
 
\bibliographystyle{amsplain}
\bibliography{subgraph_numedge}

\end{document}